\documentclass[12pt]{amsart}
\usepackage{amscd}

\def\Hn{H_n}
\def\Hnn{H_{n,n+1}}
\def\Hnp{H_{n+1}}
\def\til#1{\xrightarrow{#1}}

\let\mediumspace=\:
\let\:=\colon
\def\paragraph{\S}
\def\i{^{-1}}

\def\CC{{\mathbf C}}

\def\QQ{{\mathbf Q}}

\def\PP{{\mathbf P}}

\def\A{{\mathcal A}}

\def\F{{\mathcal F}}

\def\I{{\mathcal I}}

\def\OO{{\mathcal O}}
\def\B{{\mathcal B}}

\def\x{\times}
\def\*{\otimes}

\def\iso{\simeq}

\def\sub{\subseteq}

\def\Hilb{\operatorname{Hilb}}

\def\codim{\operatorname{codim}}

\def\tilde{\widetilde}

\swapnumbers
\newtheorem{thm}{Theorem}[section]
\newtheorem{lemma}[thm]{Lemma}
\newtheorem{prop}[thm]{Proposition}

\theoremstyle{definition}

\theoremstyle{remark}

\numberwithin{equation}{section}

\newcommand{\thmref}[1]{theorem~\ref{#1}}

\newcommand{\lemref}[1]{lemma~\ref{#1}}
\newcommand{\propref}[1]{proposition~\ref{#1}}

\begin{document}
\bibliographystyle{plain}

\title [Intersection number for punctual Hilbert sceheme]
{An intersection number for the \\ punctual Hilbert scheme of a
  surface}

\author{Geir Ellingsrud}
    \address{Mathematical Institute\\
       University of Oslo\\
       P.~O.~Box~1053\\
       N--0316 Oslo, Norway}
    \email{ellingsr@math.uio.no}
\author{Stein Arild Str{\o}mme}
     \address{Mathematical Institute\\
       University of Bergen\\
       All\'eg 55\\
       N--5007 Bergen, Norway}
     \email{stromme@mi.uib.no}

\thanks{Copyright\copyright 1996
G.~Ellingsrud and S. A. Str{\o}mme.  All rights reserved.}

\date{1996--03--21}


\maketitle

\section{Introduction}

Let $S$ be a smooth projective surface over an algebraically closed
field $k$. For any natural number $n$ let $\Hn$ denote the Hilbert
scheme parameterizing finite subschemes of $S$ of length $n$. It is
smooth and projective of dimension $2n$.

Let $P\in S$ be a point and let $M_n(P)\sub \Hn$ be the closed
reduced subvariety consisting of points which correspond to subschemes
with support at $P$. Brian\c con proved that $M_n(P)$ is an irreducible
variety of dimension $n-1$, see \cite{Bria}.

Denote by $M_n=\cup_{P\in S}M_n(P)\sub \Hn$ the subvariety whose
points correspond to subschemes with support in just one point. We may
map $M_n$ to $S$ by sending a point of $M_n$ to the point where the
corresponding subscheme is supported. The fiber of this map over a
point $P$ being the variety $M_n(P)$, we see that $M_n$ is irreducible
of dimension $n+1$.

The subvarieties $M_n$ and $M_n(P)$ are of complementary codimensions,
and hence the product of their rational equivalence classes (or dual
cohomology classes, if $k$ is the field of complex numbers) defines an
intersection number $\int_{\Hn}^{}[M_n]\cdot [M_n(P)]$. The main content
of this note is the computation of that number. The result is:

\begin{thm}\label{main}
  $\int_{\Hn}^{}[M_n]\cdot[M_n(P)]=(-1)^{n-1}n.$
\end{thm}

One reason, pointed out to us by H.~Nakajima, to be interested in these
intersection numbers is the following.  In case $k=\CC$ is the field
of complex numbers, G\"ottsche \cite{Goet1} computed the generating
series $ \sum_{m,n=0}^{\infty}\dim H^m(\Hn,\QQ)t^n u^m $ and showed
that it may be expressed in terms of classical modular forms. These
forms are closely related to the trace of some standard
representations of respectively the infinite Heisenberg algebra and
the infinite Clifford algebra.  In \cite{Naka} Nakajima defined a
representation of a product of these algebras, indexed over
$H^{*}(S,\QQ)$, on the space $\oplus_{m,n=0}^{\infty}H^{m}(\Hn,\QQ)$.
He completely described this representation up to the determination of
a series of universal constants $c_n$ for $n=1,2,\ldots$, universal in
the sense that they do not depend on the surface.  He also proved that
$c_n=\int_{\Hn}^{}[M_n]\cdot[M_n(P)]$. Hence we have

 \begin{thm}\label{NakaTheo}
   The Nakajima constants are given by $c_n=(-1)^{n-1}n$.
 \end{thm}

\section{Proof of the main theorem}

The proof will be an inductive argument comparing the number $c_n$
with $c_{n+1}$.  To make this comparison, we shall make use of the
``incidence variety'', i.e., the closed, reduced subscheme of
$\Hn\times \Hnp$ given by $\Hnn=\{(\xi,\eta)\mid \xi\sub \eta \}$.  It
is known that $\Hnn$ is smooth and irreducible of dimension $2n+2$
(there are several proofs of this, see for example \cite{Chea-1} or
\cite{Tikh-2}).

There are obvious maps $f\:\Hnn\to \Hn$ and $g\:\Hnn\to \Hnp$ induced
by the projections.  There is also a natural map $q\:\Hnn\to S$
sending a pair $(\xi,\eta)$ to the unique point where $\xi$ and $\eta$
differ.  Let $Z_n\sub \Hn\times S$ be the universal subscheme. It is
finite and flat over $\Hn$ of rank $n$.  Let $\pi_n\:Z_n\to \Hn $
denote the restriction of the projection. Furthermore, let
$\Hn'\sub\Hn$ denote the open dense subset parameterizing local
complete intersection subschemes, and put $Z_n'=\pi_n\i \Hn'$.

In the next section we will prove the following results, which
shed light on the maps $g$ and $f$.

\begin{prop}\label{structure1}
  The map $g\:\Hnn\to\Hnp$ factors naturally as
  $g=\pi_{n+1}\circ\psi$, where $\psi\:\Hnn\to Z_{n+1}$ is canonically
  isomorphic to $\PP(\omega_{Z_{n+1}})$.  In
  particular, $\psi$ is birational and an isomorphism over $Z_{n+1}'$,
  and $g$ is generically finite of degree $n+1$.
\end{prop}

\begin{prop}\label{structure2}
  The map $\phi=(f,q)\:\Hnn\to\Hn\x S$ is canonically isomorphic to
  the blowing up of $\Hn\x S$ along $Z_n$. In particular, over $Z_n'$,
  the map $\phi$ is a
  $\PP^1$-bundle.
\end{prop}

It follows that the fibers of $f$ over local complete intersection
subschemes $\xi\in\Hn'$ are given as $f\i(\xi)=\tilde S(\xi)$, the
surface $S$ blown up along $\xi$. (This is also easy to see directly.)

The locus of pairs $(\xi,\eta)\in \Hnn$ where $\xi$ and $\eta$ have
the same support is a divisor in $\Hnn$ which we denote by $E$. This
is nothing but the exceptional divisor of the blowup morphism $\phi$.
On the fiber of $f$ over a local complete intersection $\xi$, it
restricts to the exceptional divisor of $\tilde S(\xi)$.

Let $M_{n,n+1}=(g\i M_{n+1})_{\text{red}}$. We need the following
strengthening of Brian\c con's result, also to be proved in the next
section.

\begin{prop}\label{structure3}
  $M_{n,n+1}$ is irreducible, and $g$ maps it birationally to
  $M_{n+1}$. In particular, all the $M_n$ are irreducible, and the
  complete intersection subschemes form a dense open subset of $M_n$.
\end{prop}

Using the three propositions above, we have sufficient information to
carry out the intersection computation. Let us summarize the situation
in the following commutative diagram, where $B_n=\rho(M_{n,n+1})$:
\[
\begin{CD}
M_{n,n+1} @>\sub >> E    @>j>>    \Hnn     @>\gamma=(g,q)>>  \Hnp\x S \\
@VVV @V\rho VV  @VV\phi=(f,q) V @VV pr_2V \\
B_n @>\sub >> Z_n @>i>> \Hn\x S @> pr_2 >> S \\
@VVV @V\pi_n VV  @VV pr_1 V \\
M_n@>\sub>> \Hn @= \Hn
\end{CD}
\]
\begin{lemma}\label{lemmaa}
$     g^*[M_{n+1}] = (n+1)\,[M_{n,n+1}]$   in $A^n(\Hnn)$.
\end{lemma}
\begin{proof}

  Since $g\i M_{n+1}$ is a multiple structure on $M_{n,n+1}$ by
  definition, and the codimensions of $M_{n+1}$ and $M_{n,n+1}$ are
  the same, $g^*[M_{n+1}] = \ell\, [M_{n,n+1}]$ for some integer
  $\ell$. Now use that $g_*[M_{n,n+1}]=[M_{n+1}]$ (by
  \propref{structure3}) and the projection formula to get
 \[
   (n+1)[M_{n+1}] = g_* g^* [M_{n+1}] = g_* (\ell\, [M_{n,n+1}])=
   \ell\,[M_{n+1}],
 \]
 proving that $\ell=n+1$.
\end{proof}

\begin{lemma}\label{lemmab}
  $ [E]\cdot f^*[M_n] = n\,[M_{n,n+1}]$ in $A^n(\Hnn)$.
\end{lemma}
\begin{proof}
  Consider first $[M_{n,n+1}]_E\in A_{n+2}(E)$. Let $h=\pi_n\circ\rho\:
  E\to\Hn$. Since $M_{n,n+1}$ is the support of $h\i M_n$ and its
  codimension in $E$ equals $\codim(M_n,\Hn)$, we have that
  $h^*[M_n]=\ell\,[M_{n,n+1}]_E$ where $\ell$ is the multiplicity of
  $h\i M_n$ at the generic point $\eta$ of $M_{n,n+1}$. By
  \propref{structure2}, $\rho$ is a smooth at $\eta$, so $\ell$ equals
  also the multiplicity of $\pi_n\i M_n$ at the generic point
  $\rho(\eta)$ of $B_n$. But observing that $B_n$ maps isomorphically
  to $M_n$, a similar argument as in the proof of \lemref{lemmaa}
  shows that $\pi_n^*[M_n]=n\,[B_n]$, hence $\ell=n$.

  We have shown that $h^*[M_n] = n\,[M_{n,n+1}]_E$ in $A_{n+2}(E)$. Apply
  $j_*$ and the projection formula to get
  \[
  n\,[M_{n,n+1}] = j_*h^*[M_n]=j_*j^*f^*[M_n] = [E]\cdot f^*[M_n].
  \]
\end{proof}

Combining the two lemmas above, we get
\begin{equation}
  \label{vari}
  \frac1{n+1}\,g^*[M_{n+1}] = \frac1n\,[E]\cdot f^*[M_n] \in A^n(\Hnn),
\end{equation}
and exactly parallel reasoning shows that also
\begin{equation}
  \label{fast}
    \frac1{n+1}\,g^*[M_{n+1}(P)] = \frac1n\,[E]\cdot f^*[M_n(P)] \in
    A^{n+2} (\Hnn).
\end{equation}

We are now ready to prove \thmref{main}. Let $F$ be a general
fiber of $f$, for example corresponding to a reduced subscheme
$\xi$. Clearly, $f^*[M_n]\cdot f^*[M_n(P)] = c_n\,[F]$.
It is easy to see that $\int_F [E]^2 = -n$, and we get the
following computation:
\begin{alignat*}{3}
  \frac{c_{n+1}}{n+1} &= \frac{1}{n+1}\int_{\Hnp}
  [M_{n+1}][M_{n+1}(P)]&&\\[2mm]
  & =\int_{\Hnn}
  \frac1{n+1}\,{g^*[M_{n+1}]}\cdot\frac1{n+1}\,{g^*[M_{n+1}(P)]}&\quad
  &\text{(proj.\ formula)}\\[2mm]
  & = \int_{\Hnn} \frac1n\,{[E]f^*[M_n]}
  \cdot \frac1n\,{[E]f^*[M_n(P)]} &\quad&\text{(\eqref{vari} and
    \eqref{fast})}\\[2mm]
  &=c_n \int_F \frac1{n^2} \,[E]^2 =\frac{-c_n}n.
\end{alignat*}
Now since trivially $c_1=1$, \thmref{main} follows by induction.

\section{The geometry of the incidence variety}

The aim of this section is to prove propositions \ref{structure1},
\ref{structure2}, and \ref{structure3} above.  Some of the content
of this section may be found in \cite{Ell}, but for the benefit of
the reader we reproduce it here.

Consider a nested pair of subschemes $(\xi,\eta)\in\Hnn$, and let
$P=q(\xi,\eta)\in S$ be the point where they differ. There are
natural short exact sequence on $S$:
\begin{align}\label{basic0}
  &0 \to \I_{\eta} \to \I_{\xi} \to k(P) \to 0\\
  \label{basic1}
  &0 \to k(P) \to \OO_{\eta} \to \OO_{\xi} \to 0.
\end{align}
The first of these shows that the fiber $\phi\i(\xi,P)$ is
naturally identified with
the projective space $\PP(\I_{\xi}(P))$.

Dualizing \eqref{basic1} we arrive at another exact
sequence
\begin{equation}
  \label{basic2}
  0 \to \omega_{\xi} \to \omega_{\eta} \to k(P) \to 0,
\end{equation}
and this shows that the fiber $\gamma\i(\eta,P)$ maps naturally to
$\PP(\omega_{\eta}(P))$.  Dualizing again we see that \eqref{basic1}
and hence $\xi$ can be reconstructed from the right half of
\eqref{basic2}, so the map is an isomorphism.

It follows from \eqref{basic0} that
\begin{equation}
  \label{basic3}
  |\dim_k\I_{\xi}(P) - \dim_k\I_{\eta}(P)| \le 1.
\end{equation}
(If $\F$ is a coherent sheaf, $\F(P)$ means $\F\*k(P)$.)

Note also that for any pair $(\xi,P)\in
\Hn\x S$, we have
\begin{equation}
  \label{basic4}
  \dim_k \I_{\xi}(P) = 1+ \dim_k \omega_{\xi}(P),
\end{equation}
most easily seen using a minimal free resolution of the local ring
$\OO_{\xi,P}$ over $\OO_{S,P}$.

The sequences \eqref{basic0} and \eqref{basic2} can be naturally
globalized to the relative case of families of subschemes and points.
This way one easily proves \propref{structure1}, as well as the
following lemma:

 \begin{lemma}\label{philemma}
   Let $\I_n$ denote the sheaf of ideals of $Z_n$ in $\OO_{\Hn\times
     S}$.  Then there is an isomorphism $\Hnn\iso\PP(\I_n)$ such that
   $\phi$ corresponds to the tautological mapping $\PP(\I_n)\to
   \Hn\times S$.
 \end{lemma}

 We shall prove that the map $\phi\:\Hnn\to \Hn \times S$ is the blow
 up of $\Hn \times S$ along the universal subscheme $Z_n$, by proving
 a general proposition on blowing up codimension two subschemes, and
 later verify its hypotheses in the case at hand.

 Let $W$ be any irreducible algebraic scheme and $Z\sub W$ a subscheme
 of codimension $2$ whose ideal we denote by $\I_{W}$.  We assume that
 $\OO_{Z}$ is of local projective dimension 2 over $\OO_{W}$.  For any
 integer $i$ let $W_{i}=\{w\in W \mid \dim_k\I_{W}(w) = i\}$.  Let
 $\tilde W$ be the blow up of $W$ along $Z$.

 There is an obvious map from $\tilde W$ to $\PP(\I)$ due to the fact
 that $\I\OO_{\tilde W}$ is invertible.  We shall see that under
 certain conditions this map is an isomorphism. (See also
 \cite[Prop.~9]{Avramov}.)

 \begin{prop}\label{blowup}
   With the above hypothesis
      \begin{enumerate}
      \item[(a)] Suppose that $\codim W_{i}\ge i$ for all $i\ge2$. Then
        $\PP(\I)$ is irreducible and isomorphic to $\tilde W$.
      \item[(b)] If furthermore $Z$ is irreducible and $\codim
        W_{i}\ge i+1$ for $i\ge3$, then the exceptional divisor $E$ is
        irreducible.
      \end{enumerate}
 \end{prop}

\begin{proof}
The assumption on the local projective dimension gives an exact
sequence
 \begin{equation}\label{hilbertbirch}
         0\to \A \til{M} \B \to \I \to 0
 \end{equation}
 where $\A$ and $\B$ are vector bundles on $W$ whose ranks
 are $p$ and $p+1$ respectively, for some integer $p$.
 Locally the map $M$ is given
 by a $(p+1)\x p$-matrix of functions on $W$, and the
 ideal $\I$ is locally generated by its maximal minors.

 The sequence \eqref{hilbertbirch} induces a natural inclusion
 $\PP(\I)\sub\PP(\B)$.  In fact, letting $\varepsilon\:\PP(\B)\to W$ be
 the structure map and $\tau\:\varepsilon^{*}\B\to\OO_{\PP(\B)}(1)$
 the tautological quotient, $\PP(\I) $ is defined in
 $\PP(\B)$ as the vanishing locus of the map $\tau\circ M$.

 Consequently, any irreducible component of $\PP(\I)$ has dimension at
 least equal to $\dim( W)$.

 Let $\phi\:\PP(\I)\to W$ be the structure map.  Over $W_i$, the map
 $\phi$ is a $\PP^{i-1}$-bundle.  Hence by the assumption in $(a)$, we
 have
 \[ \dim \phi\i W_i \le (\dim W - i) + (i-1) < \dim W
 \]
 for all
 $i\ge 2$.  It follows that $\phi\i(W-Z)$ is dense in $\PP(\I)$,
 which is therefore irreducible.

 To see that $\PP(\I)$ is isomorphic to the blow up, we remark that it
 follows from the resolution \eqref{hilbertbirch} that
 $\I\OO_{\PP(\I)}$ is an invertible ideal.  This gives a map from
 $\PP(\I)$ to $\tilde W$, which over $\phi\i(W-Z)$ is an inverse to
 the map from $\tilde W$ to $\PP(\I)$ given above.  As both spaces are
 irreducible, the two maps are inverses to each other.

 Under the assumption in $(b)$, it follows similarly that $\phi\i
 W_2$ is dense in the exceptional locus $\phi\i Z$, as all $\phi\i W_i$
 are of strictly lower dimension if $i\ge 3$.
\end{proof}

To finish the proof of \propref{structure2},
we shall verify the conditions in proposition \ref{blowup}
for $Z=Z_{n}$ and $W=\Hn\times S$.

Let $W_{i,n}$ be the set of points $(\xi,P)\in \Hn\times S$ such that
the ideal $\I_{Z_{n}}$ needs exactly $i$ generators at $(\xi,P)$.
Equivalently,
\begin{equation}
W_{i,n}=\{(\xi,P)\in \Hn\x S \mid \dim_k\I_{\xi}(P) = i\}.
\end{equation}
We shall show by induction on $n$
that $\codim(W_{i,n},\Hn\x S) \ge 2i-2$, or
equivalently, that
\[
\dim W_{i,n}\le 2n+4-2i
\]
for all $i,n\ge1$. For $n=1$ this is evidently satisfied for all $i$.

Assume that the inequality holds for a given $n$, and all $i\ge1$.
Then it follows that
\[\dim\phi\i W_{j,n}\le (2n+4-2j) + (j-1) = 2n+3-j \le 2n+4-i\]
for all $j\ge i-1$. By \eqref{basic3},
\[\gamma\i W_{i,n+1}\sub
\phi\i W_{i-1,n}\cup\phi\i W_{i,n} \cup\phi\i W_{i+1,n}.
\]
The fibers of $\gamma$ over $W_{i,n+1}$ are $(i-2)$-dimensional. Hence
\[
\dim W_{i,n+1}+(i-2) = \dim \gamma\i W_{i,n+1}  \le 2n+4-i
\]
and hence $\dim W_{i,n+1}\le 2(n+1) +4 - 2i$, as was to be shown.

Since $2i-2\ge i$ for $i\ge2$ and $2i-2\ge i+1$ for $i\ge 3$, the
proof of \propref{structure2} is now complete. Note that we have also
proved that the exceptional divisor $E$ is irreducible.

\subsubsection*{Proof of \propref{structure3}} The only hard part is
to show that $M_{n,n+1}$ is irreducible.  We will apply
\propref{blowup} in the case where $W=M_{n}(P)\x S$ and $Z=Z_{n}\cap
W$.  As this intersection is proper, the condition on local projective
dimension still holds in this case.  Put $W_{i,n}'=W_{i,n}\cap W$.
Similar reasoning as in the last proof gives the inequality
$\codim(W_{i,n}',W)\ge i+1$ for all $n\ge 1$ and $i\ge 3$.
Now note that the exceptional divisor $\phi\i Z$ is nothing but
$M_{n,n+1}$, which is therefore irreducible by \propref{blowup}.


\end{document}